\documentclass[preprint]{imsart}
\pdfoutput=1
\RequirePackage[OT1]{fontenc}
\usepackage{amsthm,amsmath,natbib,amssymb,bm,enumerate,booktabs,float,natbib,amsxtra,enumitem}
\usepackage{yuzo_notation}
\def\mrM{{ \mathrm{M} }}
\def\mrtr{{ \mathrm{tr} }}

\RequirePackage[colorlinks,citecolor=blue,urlcolor=blue]{hyperref}
\usepackage{booktabs}

\startlocaldefs
\numberwithin{equation}{section}
\theoremstyle{plain}
\newtheorem{thm}{Theorem}[section]

\theoremstyle{remark}
\newtheorem{remark}{Remark}[section]
\newtheorem{example}{Example}[section]

\endlocaldefs
\def\citeapos#1{\citeauthor{#1}'s (\citeyear{#1})}

\allowdisplaybreaks

\begin{document}

\begin{frontmatter}
\title{An alternative to Moran's $I$ for spatial autocorrelation}
\runtitle{An alternative to Moran's $I$}

\begin{aug}
\author{\fnms{Yuzo} \snm{Maruyama}
\thanksref{t1}\ead[label=e1]{maruyama@csis.u-tokyo.ac.jp}}
\thankstext{t1}{This work was partially supported by KAKENHI \#25330035.}
\runauthor{Y. Maruyama}

\affiliation{Center for Spatial Information Science, The University of Tokyo}
\address{Center for Spatial Information Science, University of Tokyo \\
\printead{e1}}
\end{aug}
\begin{abstract}
Moran's $\II$ statistic, a popular measure of spatial autocorrelation, is revisited.
The exact range of Moran's $\II$ is given as a function of spatial weights matrix. 
We demonstrate that some spatial weights matrices lead the absolute value of upper (lower) bound larger than $1$ and that others lead the lower bound larger than $-0.5$.
Thus Moran's $\II$ is unlike Pearson's correlation coefficient.
It is also pointed out that some spatial weights matrices do not allow Moran's $\II$
to take positive values regardless of observations.
An alternative measure with exact range $[-1,1]$ is proposed
through a monotone transformation of Moran's $\II$.
\end{abstract}


\end{frontmatter}

\section{Introduction}\label{sec:intro}
Spatial autocorrelation statistics measure and analyze 
the degree of dependency among observations in a geographic space. 
In this paper, we are interested in \citeapos{Moran-1950} $\II$ statistic, which seems
the most popular measure of spatial autocorrelation.
To the best of my knowledge, however, the property has not been thoroughly investigated.
To be worse, 
there are some misunderstandings with respect to similarity to Pearson's correlation coefficient.
In some books as well as Wikipedia, it is claimed that 
the range of Moran's $\II$ is exactly $[-1,1]$ like Pearson's correlation coefficient.
However it is not true. 

In Section \ref{sec:original} of this paper, 
we give the exact range of Moran's $\II$ as a function of spatial weights matrix. 
In our numerical experiment, 
some spatial weights matrices lead the absolute value of upper (lower) bound larger than $1$. 
Others lead the lower bound larger than $-0.5$.
Thus Moran's $\II$ is quite different from Pearson's correlation coefficient
with exact range $[-1,1]$.
It is also pointed out that some spatial weights matrices do not allow Moran's $\II$
to take positive values regardless of observations.
In Section \ref{sec:modified}, we propose an alternative measure with the exact range $[-1,1]$, 
which is a monotone transformation of Moran's $\II$.

Actually, this is not the first paper which points that the range of Moran's $\II$ is not $[-1,1]$. 
To the best of my knowledge, \cite{Cliff-Ord-1981}, \cite{Goodchild-1988}, \cite{deJong-etal-1984} and
\cite{Waller-Gotway-2004} did so. 
From mathematical viewpoint, Theorem \ref{thm:moran_original} 
in Section \ref{sec:original} is essentially the same as the results by \cite{deJong-etal-1984}.
However any alternative with exact range $[-1,1]$, which I hope, 
helps better understanding of spatial autocorrelation, has not yet been proposed. 
Furthermore it is notable that since the alternative is proposed by 
a monotone transformation of Moran's $\II$, 
the permutation test of spatial independence based on Moran's $\II$, 
which has been implemented in earlier studies, 
can be saved after the standard measure is replaced by the alternative.

\section{Moran's $\II$}\label{sec:original}

In this section, we give the exact range of Moran's $\II$.
Suppose that, in $n$ spatial units $s_1,\dots,s_n$, we observe $y_1=y(s_1), \dots, y_n=y(s_n)$.
We also suppose that spatial weights $w_{ij}$, the weights 
between each pair of spatial units $s_i$ and $s_j$, 
which satisfy
\begin{align}\label{def_w}
w_{ij}\geq 0\text{ for any }i\neq j,\quad w_{ii}=0 \text{ for any }i,
\end{align}
are given by some preset rules. 
Then Moran's $\II$ statistic is given by
\begin{align*}
I=\frac{n\sum_{i=1}^n \sum_{j=1}^n w_{ij}(y_i-\bar{y})(y_j-\bar{y})}
{\{\sum_{i,j=1}^n w_{ij}\}\sum_{i=1}^n(y_i-\bar{y})^2}. 
\end{align*}
Let $\bmy=(y_1,\dots,y_n)^\T$ and
$\tilde{\bmy}=\bmy-\bar{y}\bmone_n=(\bm{I}_n-\bmone_n\bmone_n^\T/n)\bmy$
where $^\T$ denotes the transpose and $\bmone_n$ denotes the $n$-dimensional vector of ones.
With the spatial weights matrix $\bmW=(w_{ij})$, 
Moran's $\II$ statistic is rewritten as
\begin{align*}
\II=\frac{n}{\sum_{i,j=1}^n w_{ij}}
\frac{\tilde{\bmy}^\T\bmW\tilde{\bmy}}{\tilde{\bmy}^\T\tilde{\bmy}}.
\end{align*}
When $\bmW$ is not symmetric, 
the symmetric alternative $(\bmW+\bmW^\T)/2$ can be identified to $\bmW$ since
\begin{gather*}
 \tilde{\bmy}^\T\bmW\tilde{\bmy}=\tilde{\bmy}^\T\{(\bmW+\bmW^\T)/2\}\tilde{\bmy}, 
\text{ for any }\tbmy, \\
\sum\nolimits_{i,j=1}^n w_{ij}=\bmone_n^\T\bmW\bmone_n=
\bmone_n^\T\{(\bmW+\bmW^\T)/2\}\bmone_n=\sum\nolimits_{i,j=1}^n \{w_{ij}+w_{ji}\}/2.
\end{gather*}
Thus, in this paper, we assume that $\bmW$ is symmetric.

Let $\mathcal{V}$ be the $(n-1)$-dimensional subspace of $\mbR^n$ orthogonal to $\bmone_n$.
Then $\bm{I}_n-\bmone_n\bmone_n^\T/n$ is the matrix for projection onto $\mathcal{V}$.
The spectral decomposition of $\bm{I}_n-\bmone_n\bmone_n^\T/n$ is given by
\begin{align}\label{eq:HHT}
 \bm{I}_n-\bmone_n\bmone_n^\T/n=
\sum\nolimits_{i=1}^{n-1}\bmh_i\bmh_i^\T=\bmH\bmH^\T
\end{align}
where $\bmh_i\in\mbR^n$, $\bmh_i^\T\bmh_i=1$ and $\bmh_i^\T\bmone_n=0$ for any $1\leq i\leq n-1$.
Further, in \eqref{eq:HHT}, 
$ \bmH=(\bmh_1,\dots,\bmh_{n-1})$ is the $n\times(n-1)$ matrix where
$\{\bmh_1,\dots,\bmh_{n-1}\}$ span the subspace $\mathcal{V}$ as the orthonormal bases.

Let $\bmv=\bmH^\T \bmy\in\mbR^{n-1}$. 
Then Moran's $\II$ statistic is rewritten as
\begin{equation}\label{eq:new_QF}
\begin{split}
\II= \frac{\tilde{\bmy}^\T\bmW\tilde{\bmy}}{n\bar{w}\tilde{\bmy}^\T\tilde{\bmy}}
=\frac{\bmy^\T\bmH\bmH^\T\bmW\bmH\bmH^\T\bmy}{n\bar{w}\bmy^\T\bmH\bmH^\T\bmy}
=\frac{\bmv^\T\tbmW\bmv}{\bmv^\T\bmv}
\end{split}
\end{equation}
where $\bar{w}=\sum_{i,j=1}^n w_{ij}/n^2$ and $\tbmW$ is the $(n-1)\times (n-1)$ matrix given by
\begin{equation}\label{eq:tbmW}
 \tbmW=(n\bar{w})^{-1}\bmH^\T\bmW\bmH.
\end{equation}
For \eqref{eq:new_QF}, the standard result of maximization and minimization of quadratic form
is available. 
Let $\lambda_{(1)}$ and $\lambda_{(n-1)}$ be the smallest and largest eigenvalue, 
among $n-1$ eigenvalues of $\tbmW$, respectively. Then, for any $\bmv\in\mbR^{n-1}$,
\begin{align*}
\II\in\left[\lambda_{(1)},\lambda_{(n-1)}\right].
\end{align*}
The lower and upper bound of $\II$ is attained by
\begin{align*}
\II=\begin{cases}
 \lambda_{(1)} & \bmv\propto \bmj_{(1)} \\ 
 \lambda_{(n-1)} & \bmv\propto \bmj_{(n-1)},
 \end{cases}
\end{align*}
where $\bmj_{(1)}$ and $\bmj_{(n-1)}$ are the eigenvectors corresponding to 
$\lambda_{(1)}$ and $\lambda_{(n-1)}$ respectively.
By the following equivalence,
\begin{align*}
\bmv\propto \bmj \ &\Leftrightarrow \ \bmH^\T\bmy \propto \bmj \ \Leftrightarrow \ \bmH\bmH^\T\bmy\propto \bmH\bmj \ \Leftrightarrow \ (\bmI-\bmone_n\bmone_n^\T/n)\bmy \propto \bmH\bmj \\
&\Leftrightarrow \ \bmy=c_1\bmone_n+c_2\bmH\bmj \text{ for any }c_1 \text{ and nonzero } c_2,
\end{align*}
we have a following result.

\begin{thm}\label{thm:moran_original}
The exact range of Moran's $\II$ is 
\begin{equation*}
\II\in\left[\lambda_{(1)},\lambda_{(n-1)}\right]
\end{equation*}
where 
$\lambda_{(1)}$ and $\lambda_{(n-1)}$ be the smallest and largest eigenvalue of
among $n-1$ eigenvalues of $\tbmW$ given by \eqref{eq:tbmW}.
The lower and upper bound is attained as follows;
\begin{equation*}
\II= \begin{cases}
 \lambda_{(1)} & \bmy=c_1 \bmone_n +c_2 \bmH\bmj_{(1)}\\
 \lambda_{(n-1)} & \bmy=c_3 \bmone_n +c_4 \bmH\bmj_{(n-1)},
 \end{cases}
\end{equation*}
where any $c_1,c_3$ and any $c_2,c_4\neq 0$.
\end{thm}
\begin{remark}
In \eqref{eq:HHT}, the choice of $\bmH$ or $\bmh_1,\dots,\bmh_{n-1}$ is arbitrary.
Clearly $ \lambda_{(1)}$ and $\lambda_{(n-1)}$ are also expressed by
\begin{align*}
 \lambda_{(1)}=\min_{\bmz^\T\bmone_n=0}\frac{\bmz^\T\bmW\bmz}{n\bar{w}\bmz^\T\bmz},\quad
 \lambda_{(n-1)}=\max_{\bmz^\T\bmone_n=0}\frac{\bmz^\T\bmW\bmz}{n\bar{w}\bmz^\T\bmz}
\end{align*}
which do depend on $\bmW$ not on $\bmH$. In practice, when $\lambda_{(1)}$ and $\lambda_{(n-1)}$
are calculated, the choice of $\bmH=(\bmh_1,\dots,\bmh_{n-1})$ 
based on Helmert orthogonal matrix is available, that is,
$\bmh_i$ for $i=1,\dots,n-1$ given by
\begin{align*}
 \bmh_i^\T=(\bmone_{i}^\T,-i,\bmzero_{n-i-1}^\T)/\sqrt{i(i+1)}.
\end{align*}
\end{remark}

\begin{example}
For simplicity, suppose $s_1,\dots,s_n$ are on a line with equal spacing
and spatial weights for $s_i,s_j$ is
\begin{equation*}
w_{ij}=
\begin{cases}
 2^{-|i-j|+1} & 1\leq |i-j|\leq q \\
 0 & i=j, \  |i-j|> q.
\end{cases}
\end{equation*}
Table \ref{table_1} gives lower and upper bound of Moran's $\II$ for some $n$ and $q$.
For any $n$ and $q$, the range of $\II$ is not $[-1,1]$.
It is demonstrated that 
some spatial weights matrices lead the absolute value of upper (lower) bound larger than $1$
and that others lead the lower bound larger than $-0.5$, which is unexpected.
Thus Moran's $\II$ is unlike Pearson's correlation coefficient with exact range $[-1,1]$. 
\end{example}

\begin{table}
\caption{attainable lower \& upper bound of Moran's $\II$}
\centering
\begin{tabular}{ccccccccc} \toprule
$n\backslash q$ & \multicolumn{2}{c}{1} & & \multicolumn{2}{c}{2} & & \multicolumn{2}{c}{3} \\ \midrule
& lower & upper & & lower & upper & & lower & upper \\ 
10 & -1.066 & 0.935 & & -0.541 & 0.831 & & -0.482 & 0.746 \\
20 & -1.041 & 1.006 & & -0.526 & 0.981 & & -0.457 & 0.955 \\
30 & -1.029 & 1.013 & & -0.519 & 1.005 & & -0.449 & 0.995 \\
40 & -1.023 & 1.014 & & -0.514 & 1.011 & & -0.444 & 1.006 \\
50 & -1.018 & 1.013 & & -0.512 & 1.012 & & -0.441 & 1.010 \\ \bottomrule
\end{tabular}
\label{table_1}
\end{table}
\begin{remark}
In all cases in Table \ref{table_1}, both $\lambda_{(1)}<0$ and $\lambda_{(n-1)}>0$ are satisfied.
Researchers and practitioners would like measure of spatial correlation 
to take both positive and negative values depending upon positive and negative
spatial correlation. 
It is clear that indefiniteness of $\tbmW$ is equivalent to $\lambda_{(1)}<0$ and $\lambda_{(n-1)}>0$ simultaneously.
Actually, as in below, some spatial weights matrices lead $\lambda_{(n-1)}<0$, 
or negative definiteness of $\tbmW$.
In such a case, Moran's $\II$ cannot take positive values for any $\bmy$, which is not desirable
at all.

In order to investigate definiteness of $\tbmW$, 
the following properties are useful on eigenvalues, trace or sum of all diagonal elements and
their relationship;
\begin{enumerate}[font=\bfseries, label=P\arabic*]
\item \label{trace_eigen_2}
When $\bm{AB}$ and $\bm{BA}$ are square matrices, $\mrtr(\bm{AB})=\mrtr(\bm{BA})$,
\item \label{trace_eigen_3}
When both $\bm{A}$ and $\bm{B}$ are square matrices,
$\mrtr(\bm{A}+\bm{B})=\mrtr\bm{A}+\mrtr\bm{B}$,
\item \label{trace_eigen_1}
Trace is equal to sum of all eigenvalues.
\end{enumerate}
From \ref{trace_eigen_2}, \ref{trace_eigen_3} and $\sum w_{ii}=0$, we have
\begin{equation}\label{tr_tbmW}
\begin{split}
 \mrtr\tbmW
&=(n\bar{w})^{-1}\mrtr\bmH^\T\bmW\bmH =(n\bar{w})^{-1}\mrtr\bmW\bmH\bmH^\T \\
&=(n\bar{w})^{-1}\mrtr\bmW(\bm{I}-\bmone_n\bmone_n^\T/n) =(n\bar{w})^{-1}\mrtr\bmW - (n^2\bar{w})^{-1}\mrtr\bmone_n^\T\bmW\bmone_n \\
& =(n\bar{w})^{-1}\sum\nolimits_{i=1}^n w_{ii} - (n^2\bar{w})^{-1}\sum\nolimits_{i,j=1}^nw_{ij} =-1.
\end{split}
\end{equation}
From \ref{trace_eigen_1} together with $\mrtr\tbmW=-1$, 
$\lambda_{(1)}<0$ are guaranteed whereas $\lambda_{(n-1)}>0$ are not.

In fact, if $w_{ij}=1$ for any $i\neq j$ or $\bmW=\bmone_n\bmone_n^\T-\bm{I}$, we have
$n\bar{w}=n-1$ and
\begin{align*}
 \tbmW&=(n\bar{w})^{-1}\bmH^\T\bmW\bmH=(n-1)^{-1}\bmH^\T(\bmone_n\bmone_n^\T-\bm{I})\bmH \\
&=-(n-1)^{-1}\bmH^\T\bmH=-(n-1)^{-1}\bmI_{n-1}
\end{align*}
whose eigenvalues are all equal to $-(n-1)^{-1}$.
This implies that $ \tbmW$ is negative definite and that Moran's $\II$ is negative for any $\bmy$.

The negative definiteness with $w_{ij}=1$ for any $i\neq j$
suggests possible negative definiteness for $w_{ij}$ with small variability.
Suppose, for $0<a<1$, we observe $n(n-1)$ uniform random variables on $(1-a,1+a)$,
which are assigned to $w_{ij}$ for $i\neq j$.
We are interested in definiteness of 
\begin{equation*}
\tbmW= (n\bar{w})^{-1}\bmH^\T\left\{(\bmW+\bmW^\T)/2\right\}\bmH.
\end{equation*}
From numerical experiment, we have negative definiteness and indefiniteness of $\tbmW$ for $a\in (0,a_*)$
and $a\in(a_*,1)$, respectively, where
$a_*$ is approximately equal to 
$0.3$ for $n=25$; $0.2$ for $n=50$; $0.14$ for $n=75$; $0.12$ for $n=100$.
\end{remark}

\section{An alternative to Moran's $\II$}\label{sec:modified}
In section \ref{sec:original}, we pointed out some disadvantages of the original Moran's $\II$.
In this section, we propose an alternative to Moran's $\II$, the range of which is exactly $[-1,1]$
and which is a monotone transformation of the original Moran's $\II$.
There are two steps of modification. The one is a linear transformation 
for acquirement of indefiniteness.
The other is an adequate normalization for the linear transformation.
 
As pointed out in \eqref{tr_tbmW}, the sum of eigenvalues of $\tbmW=(n\bar{w})^{-1}\bmH^\T\bmW\bmH$
is $-1$, which seems related to the bias of the original Moran's $\II$,
that is, 
\begin{align*}
 \EE[\II]=-\frac{1}{n-1}<0
\end{align*}
under spatial independence, as in \cite{Cliff-Ord-1981}.
In order to get a kind of unbiasedness,
we start with a linear transformation of Moran's $\II$,
\begin{align*}
(n-1)\II+1
=(n-1)\frac{\bmv^\T\tbmW\bmv}{\bmv^\T\bmv}+1=\frac{\bmv^\T\bmQ\bmv}{\bmv^\T\bmv},
\end{align*}
where $\bmv=\bmH^\T\bmy\in\mbR^{n-1}$ and $\bmQ=(n-1)\tbmW+\bmI_{n-1}$.
Since we have
\begin{align*}
 \mrtr\bmQ=(n-1)\mrtr\tbmW+\mrtr\bmI_{n-1}=-(n-1)+(n-1)=0,
\end{align*}
the $(n-1)\times (n-1)$ matrix $\bmQ$ is indefinite, that is,
$(n-1)\II+1$ can take both positive and negative values depending on $\bmy$.

We are now in position to give the adequate normalization of $(n-1)\II+1$.
The smallest and largest eigenvalue of $\bmQ$ is $(n-1)\lambda_{(1)}+1$
and $(n-1)\lambda_{(n-1)}+1$, which are guaranteed to be negative and positive respectively,
by the indefiniteness of $\bmQ$ or $\mrtr\bmQ=0$.
Hence, 
for $\bmv$ which satisfies $\bmv^\T\bmQ\bmv<0$ or equivalently $(n-1)\II+1<0$, we have
\begin{equation*}
\frac{-\bmv^\T\bmQ\bmv}{\bmv^\T\bmv}\leq \left|(n-1)\lambda_{(1)}+1\right|
\end{equation*} 
and the equality is attained by $\II=\lambda_{(1)}$.
For $\bmv$ which satisfies $\bmv^\T\bmQ\bmv>0$ or equivalently $(n-1)\II+1>0$, we have
\begin{equation*}
\frac{\bmv^\T\bmQ\bmv}{\bmv^\T\bmv}\leq (n-1)\lambda_{(n-1)}+1
\end{equation*} 
and the equality is attained by $\II=\lambda_{(n-1)}$.

Let $\II_{\mrM}$ be an alternative to Moran's $\II$ given by
\begin{align*}
 \II_{\mrM}=\frac{(n-1)\II+1}{C}, \quad C=
\begin{cases}
|(n-1)\lambda_{(1)}+1| & (n-1)\II+1<0 \\
(n-1)\lambda_{(n-1)}+1 & (n-1)\II+1\geq 0 
\end{cases}
\end{align*}
where $ \lambda_{(1)}$ and $ \lambda_{(n-1)}$ are the smallest and largest eigenvalue
of the $(n-1)\times (n-1)$ matrix $\tbmW$ given by \eqref{eq:tbmW}.
Then we have a following result.
\begin{thm}\label{thm:moran_modified}
The exact range of $\II_{\mrM}$ is 
\begin{equation*}
\II_{\mrM}\in\left[-1,1\right].
\end{equation*}
The lower and upper bound is attained as follows;
\begin{equation*}
\II_{\mrM}= \begin{cases}
 -1 & \II=\lambda_{(1)} \\
 1 & \II=\lambda_{(n-1)}. 
 \end{cases}
\end{equation*}
\end{thm}
As far as we know, this is the first measure with exact range $[-1,1]$, which I hope, 
helps better understanding of spatial autocorrelation.
Furthermore it is notable that since the alternative $\II_{\mrM}$ is
a monotone transformation of Moran's $\II$, 
the permutation test of spatial independence based on Moran's $\II$, 
which has been implemented in earlier studies, 
can be saved after the standard measure is replaced by $\II_{\mrM}$.

\end{document}